\documentclass[journal]{IEEEtran}

\usepackage{amssymb,amsmath,amsthm,graphicx,verbatim,amsbsy,euscript}
\usepackage{enumerate}
\usepackage{mathtools}
\usepackage{fullpage}
\usepackage{graphicx}
\usepackage[center]{caption}
\usepackage{cancel}
\usepackage{subfig}
\usepackage{epsfig} 
\usepackage{epstopdf}
\usepackage{multirow}

\usepackage{graphicx}      

\usepackage{graphics} 
\usepackage{epsfig} 
\usepackage{epstopdf}
\usepackage{subfig}

\usepackage{caption}
\usepackage{mathptmx} 
\usepackage{times} 

\usepackage{amsmath} 
\usepackage{amssymb}  
\usepackage{algorithm}

\usepackage{algpseudocode}

\newtheorem{theorem}{Theorem}[section]
\newtheorem{lemma}[theorem]{Lemma}

\theoremstyle{definition}

\theoremstyle{remark}

\numberwithin{equation}{section}

\newcommand{\p}{\ensuremath{\partial}}

\newcommand{\R}{\mathbb{R}}

\newcommand{\f}[2]{\ensuremath{\frac{#1}{#2}}}

\newcommand{\norm}[1]{\left\|#1\right\|}
\newcommand{\abs}[1]{\left|#1\right|}
\newcommand{\br}[1]{\left(#1\right)}

\usepackage{cite}

\ifCLASSINFOpdf
  
\else
  
\fi

\begin{document}

\title{Multi-Robot Control Using Time-Varying Density Functions}

\author{Sung~G.~Lee
        and~Magnus~Egerstedt,~\IEEEmembership{Fellow,~IEEE}
\thanks{The authors are with the Department
of Electrical and Computer Engineering, Georgia Institute of Technology, Atlanta,
GA, 30332 USA (e-mail: slee656@gatech.edu; magnus@gatech.edu}
\thanks{Manuscript received ?; revised ?.}}

\markboth{IEEE Transactions on robotics, vol. ?, no.?, ? ?}%
{Shell \MakeLowercase{\textit{et al.}}: Bare Demo of IEEEtran.cls for Journals}

\maketitle

\begin{abstract}
This paper presents an approach to externally influencing a team of robots by means of time-varying density functions. These density functions represent rough references for where the robots should be located. To this end, a continuous-time algorithm is proposed that moves the robots so as to provide optimal coverage given the density functions as they evolve over time. The developed algorithm represents an extension to previous coverage algorithms in that time-varying densities are explicitly taken into account in a provable manner.  A distributed approximation to this algorithm is moreover proposed whereby the robots only need to access information from adjacent robots.  Simulations and robotic experiments show that the proposed algorithms do indeed exhibit the desired behaviors in practice as well as in theory.
\end{abstract}


\IEEEpeerreviewmaketitle

\section{Introduction}

In this paper, we present a novel approach to influencing a team of robots using optimal coverage algorithms for general time-varying density functions. This has potential implications for how human operators can interact with large teams of mobile robots, where one of the main challenges is the construction of suitable abstractions that make the entire team amenable to human control, e.g., \cite{Kolling:HSI}. For such abstractions to be useful, they need to scale gracefully as the number of robots increases. As such, density functions are promising such abstractions in that they are independent of the team size. 

Coverage control is one area of multi-agent control that has recieved significant attention lately, e.g., \cite{Martinez:Motion}, \cite{Ghosh:Coverage},  and it is concerned with how to position agents in such a way that ``surveillance" of a domain of interest is maximized. In this context, an idea that has been 
widely adopted to describe how interesting a "domain of interest" is, is to associate a density function to the domain, as was done in \cite{Cortes:Coverage, Cortes:Coverage2, Cortes:Coordination, Pimenta:SCAT, Haumann:Discoverage, Martinez:Motion}. However, the focus of previous coverage algorithms
has  largely been on static density functions, which does not provide enough flexibility when human operators are to adaptively interact with a team through a dynamic re-shaping of the density functions.

To enable this line of inquiry, we require an algorithm that can guarantee multi-robot optimal coverage given general time-varying density functions. Applications to this beyond the means for multi-robot influence can be found in a number of domains. For example, in  search and rescue scenarios, the density function could represent the probability of a lost person being at a certain point in an area, e.g., \cite{Macwan:Target}. Additionally,  optimal coverage of density functions for multi-robot surveillance and exploration was used in \cite{Haumann:Discoverage}, where the density function was modeled to be a function of the explored "frontier."  (For other examples, see \cite{Ghaffarkhah:DynamicCoverage} and references therein.)

To date, 
relatively little work has been done on coverage with time-varying density functions. In \cite{Cortes:Coverage}, the time-varying case was investigated under a set of simplifying assumptions on the density functions. While the resulting algorithm
works well for many choices of density functions, we will show in later sections that the assumptions in \cite{Cortes:Coverage} do not hold in general and  their violations may even cause the algorithm to break down. Another stab at the problem was pursued in \cite{Pimenta:SCAT}, where time-varying density functions where used as a means to tracking moving targets. While simulations and experiments verified that coverage was indeed achieved, formal guarantees were absent.

In contrast to \cite{Cortes:Coverage} and \cite{Pimenta:SCAT}, in this paper we derive an algorithm that guarantees optimal coverage with quite general, time-varying functions. As with the algorithms for time-invariant density functions, CVT (centroidal Voronoi tessellations) will play a key role.
A CVT is a configuration where the positions of each robot coincide with the centroids of their Voronoi cells, given a so-called Voronoi tessellation of the domain of interest. The algorithm proposed in this paper will achieve optimal coverage by first letting the robots converge to a CVT associated with a static density function as an initialization step, and then have the robots maintain the CVT as the density function starts evolving over time.

This paper is organized as follows: In Section \ref{sec:loc} the problem setup is discussed in the context of locational costs that evaluate how effective given robot configurations are at achieving coverage. This is followed by the formulation of the main, centralized algorithm for coverage with time-varying density functions in Section \ref{sec:time}. This centralized algorithm is approximated in a decentralized manner using a truncated Neumann series in Section \ref{sec:Decentralization}. The different algorithms are implemented and compared on five mobile robots in Section \ref{sec:rob}.


\section{LOCATIONAL COSTS AND VORONOI TESSELLATIONS}
\label{sec:loc}

In order to even talk about optimal coverage, one first has to associate a cost to a robot configuration that describes how well a given area is being covered. For this, we will follow the construction of this so-called locational cost, as was done, for example, in \cite{Cortes:Coverage}, and we stress that no results in this section are new -- we simply include them for the sake of easy reference.  

Let $D \subset \R^2$ be the two-dimensional convex domain representing the area of interest. Moreover, let $\phi:D\times [0, \infty) \rightarrow (0, \infty)$ be the associated density function, which we will assume is bounded and continuously differentiable in both arguments, and where $\phi(q,t)$ captures the relative importance of a point $q\in D$ at time $t$.

Now, the coverage problem involves placing $n$ robots in $D$, and we let $p_i \in D$, $i=1, \cdots, n$ be the position of the $i$th robot. Moreover, the domain itself will be divided into regions of dominance, e.g., \cite{Cortes:Coverage2}, $P_1,\ldots,P_n$ (forming a proper partition of $D$), where the idea is to let robot $i$ be in charge of covering region $P_i$. One can then ask how good the choice of $p$ and $P$ is, where $p=[p_1^T,\ldots,p_n^T]^T$, and $P=\{P_1, \cdots, P_n\}$. The final piece needed to answer this question is a measure of how well a given point $q\in D$ is covered by robot $i$ at position $p_i\in D$ (see \cite{Du:CVT} and references therein). 
As the
performance of a large class of sensors deteriorate with a rate proportional to the square of the distance \cite{Meguerdichian:Exposure}, \cite{Adlakha:Critical}, the resulting locational cost is
\begin{equation}
H(p,P,t) = \sum_{i=1}^{n} \int_{P_i} \norm{q-p_i}^2\phi(q,t)dq. \label{cost not voronoi}
\end{equation}

At a given time $t$, when a configuration of robots ($p$) together with the partition ($P$) minimize \eqref{cost not voronoi}, the domain is said to be optimally covered
 with respect to $\phi$. However, it is possible to view the minimization problem as a function of $p$ alone, \cite{Cortes:Coverage},
by observing that
%
%
given $p$, the choices of $P_i$ that minimize \eqref{cost not voronoi} is 
\[
V_i(p) = \{q \in D \,\, | \,\, \norm{q-p_i} \leq \norm{q-p_j}, i \neq j \}.
\]
This partition of $D$ is a Voronoi tessellation -- hence the use of $V_i$ to denote the region. With this choice of region, we can remove the partition as a decision variable and instead focus on the locational cost
\begin{equation}
H(p,t) = \sum_{i=1}^{n} \int_{V_i(p)} \norm{q-p_i}^2 \phi(q,t)dq  \label{cost}
\end{equation}

In \cite{Iri:Fast,Du:CVT} it was shown that
\begin{equation}
\frac{\partial H}{\partial p_i} = \int_{V_i} -2(q - p_i)^T \phi(q,t)dq,
\label{eq:part}
\end{equation}
and since $\phi > 0$, one can define the mass $m_i$ and center of mass $c_i$ of the $i$-th Voronoi cell, $V_i$, as
\begin{equation}
m_i (p,t) = \int_{V_i(p)} \phi(q,t)dq \label{M def}
\end{equation}
\begin{equation}
c_i (p,t) = \frac{ \int_{V_i(p)} q \phi(q,t)dq}{m_i} \label{CM def}.
\end{equation}
Using these quantities, the partial derivative in Equation \ref{eq:part} can be rewritten as
\begin{equation}
\frac{\partial H}{\partial p_i} = 2m_i (p_i - c_i)^T.
\label{eq:partH}
\end{equation}
From this expression, we can see that a critical point of \eqref{cost} is
\begin{equation}
p_i(t) = c_i(p,t), \quad i=1,\cdots,n,
\label{eq:CVT}
\end{equation}
and a minimizer to \eqref{cost} is necessarily of this form, \cite{Du:Convergence}. Moreover, when Equation \ref{eq:CVT} is satisfied, $p$ is a so-called centroidal Voronoi tessellation (CVT). 

The robots being in a CVT configuration does not, however, imply that the global minimum of \eqref{cost} is attained. In fact, the CVT is in general not unique given a density function $\phi$.\begin{footnote}{\cite{Du:CVT} gives an example where two robots can be in multiple different CVTs with different coverage costs, with respect to the same density function.}\end{footnote} Finding the globally minimizing configuration is a difficult problem due to the nonlinearity and nonconvexity of \eqref{cost}, as discussed in \cite{Liu:CVT}. As such, in this paper, we are interested in designing algorithms that guarantee convergence to local minima with respect to time-varying density functions, and we make no claims about finding the global minimum.

In light of Equation \ref{eq:partH}, the gradient direction (with respect to $p_i$)  is given by $(p_i-c_i)$. As such, a (scaled) gradient descent motion for the individual robots to execute would be

{\vspace{1ex} \noindent \fbox{ \begin{minipage}{.95\columnwidth} 
\noindent
{\it Lloyd:}
\begin{equation}
\dot p_i = -k(p_i - c_i)
\label{eq:Lloyd}
\end{equation}
\end{minipage}}\\[2ex]}
where $k$ is a positive gain.
This is a continuous-time version of Lloyd's algorithm for obtaining CVTs as long as $\phi$ does not depend on $t$.
The way to see this, as was done in \cite{Cortes:Coverage2}, is to take $H(p)$ in  Equation \eqref{cost} (note that we assume that $H$ only depends on $p$ and not on $t$ for the purpose of this argument) as the Lyapunov function,
\begin{align*}
\frac{d}{dt} H(p) & = \sum_{i=1}^{n}\frac{\partial}{\partial p_i} H(p) \dot{p_i}\\
& = \sum_{i=1}^{n}2m_i (p_i - c_i)^T (-k(p_i - c_i))\\
& = -2k\sum_{i=1}^{n}m_i \norm{p_i - c_i}^2
\end{align*}
By LaSalle's invariance principle, the multi-robot system asymptotically converges to a configuration $\{\norm{p_i - c_i}^2=0,\,\, i=1,\cdots,n \}$, i.e., to a CVT, \cite{Cortes:Coverage2}. 

However, if $\phi$ is time-varying, the same control law does not stabilize the multi-robot system to a CVT. This point can be hinted at by investigating the evolution of a time-dependent $H(p,t)$,
\begin{align*}
\frac{d}{dt} & H(p,t)  = \sum_{i=1}^{n}\frac{\partial}{\partial p_i} H(p,t) \dot{p_i} + \frac{\partial}{\partial t} H(p,t)\\
 = & \sum_{i=1}^{n} \int_{V_i}  \norm{q-p_i}^2 \frac{\partial \phi}{\partial t}(q,t)dq  -2k\sum_{i=1}^{n}m_i \norm{p_i - c_i}^2.
\end{align*}
There is no reason, in general, to assume that this expression is negative since we do not want to impose assumptions on slowly varying density functions. Instead, what is needed is a new set of algorithms for handling the time-varying case, which is the topic of the next section.

\section{TIME-VARYING DENSITY FUNCTIONS}
\label{sec:time}

To get around the problem associated with non-slowly varying density functions, timing information must be included in the motion of the robots. In \cite{Cortes:Coverage}, this was done through the assumption that $\phi(q,t)$ is such that 
	\[
	\f{d}{dt} \br{ \sum_{i = 1}^{n} \int_{V_i}\norm{q-c_i}^2\phi(q,t)dq}=0.
\]
Letting
\[
m_{i,t} = \int_{V_i} \dot \phi(q,t)dq, \quad c_{i,t} = \f{1}{m_i}\br{\int_{V_i}q\dot\phi(q,t)dq- m_{i,t} c_i},
\]
the algorithm in \cite{Cortes:Coverage} for time-varying density functions is given by

{\vspace{1ex} \noindent \fbox{ \begin{minipage}{.95\columnwidth} 
\noindent
{\it Cortes:}
\begin{equation}
\dot p_i = c_{i,t} -(k + \f{m_{i,t}}{m_i})(p_i - c_i).
\label{eq:Cortes}
\end{equation}
\end{minipage}}\\[2ex]}
Under the  previously mentioned assumption on $\phi$, $H(p,t)$ again becomes a Lyapunov function when the agents move according to Equation \ref{eq:Cortes}, and convergence to a time-varying CVT is established. 

Unfortunately, the assumption required to make Equation \ref{eq:Cortes} work is rather restrictive and for the remainder of the paper, we will develop new methods for handling time-varying density functions without having to make restrictive assumptions. The reason why we do not want to impose these assumptions on $\phi(q,t)$ is that the density function is to be thought of as an external, human-generated input to the system. And there are no a priori reasons why the human operator would restrict the interactions to satisfy particular  regularity assumptions on $\phi$.


One way forward is to note that if we are already at a CVT at time $t_0$, i.e., $p(t_0)=c(p(t_0),t_0))$, where $c=[c_1^T,\ldots,c_n^T]^T$, it should be possible to maintain the CVT. In other words, if we can enforce that 
\[
\f{d}{dt} ( p(t)-c(p(t),t) ) = 0 \quad \forall t \geq t_0,
\]
the time-varying CVT would have been maintained. This means that
$$
\dot p=\dot c=\frac{\partial c}{\partial p}\dot p+\frac{\partial c}{\partial t},
$$
which rearranges to
\begin{equation}
\dot p = \left( I-\frac{\partial c }{ \partial p } \right)^{-1}\frac{\partial c }{ \partial t }.
\label{eq:inv}
\end{equation}
As such, we have established the following result

\begin{theorem}
\textit{Let $p(t_0) = c(p(t_0),t_0)$. If  
\[
\dot p = \left( I-\frac{\partial c }{ \partial p } \right)^{-1}\frac{\partial c }{ \partial t },~t\geq t_0
\]
then
\[
\norm{p(t)-c(p(t),t)} = 0, ~ t \geq t_0
\]
as long as the inverse $(I-\partial c/\partial p)^{-1}$ is well-defined.}
\label{thm:thm1}
\end{theorem}

There are a number of issues that must be resolved about the evolution in Equation \ref{eq:inv}, namely $(i)$ When is the inverse well-defined?; 
$(ii)$ How can one ensure that $p(t_0)=c(p(t_0),t_0)$?; $(iii)$ How is $\partial c/\partial p$ computed?; and $(iv)$ Is it possible to implement this in a distributed manner? The  first question is in general quite hard to answer. In \cite{Du:Acceleration} it was shown that in the time-invariant case, the inverse is well-defined as long as $\phi(p)$ is a log-concave function of $p$. Moreover, we need $\phi$ continuously differentiable in  both arguments, so these two conditions are enough to ensure that the inverse exists. However, this is not particularly satisfying and it does indeed pose a major challenge to the ambition of providing algorithms that can handle general, time-varying density functions. As will be seen in Section \ref{sec:Decentralization}, it is possible to get around this restriction while, at the same time, answer the fourth question through the introduction of a well-posed Neumann approximation of the inverse as a  mechanism for achieving distributed versions of the algorithm. The answer to the remaining two questions will be discussed below.

%
%


The first issue to be addressed is the constraint that $p(t_0)=c(p(t_0),t_0)$ for some initial time $t_0$. This is, practically speaking, easily achievable
since we treat the density function $\phi$ as an input to the multi-robot system. As such, we can pick static $\phi$  initially and deploy a time-invariant algorithms, such as Lloyd's algorithm, until the robots achieve a CVT.\begin{footnote}{This really happens asymptotically, but practically it is enough to get sufficiently close to a CVT.}\end{footnote}  Once the robots are at a CVT, we can allow $\phi$ to vary with time and employ the motion in Equation \ref{eq:inv}. Or, even better, add a proportional term that forces the robots to a CVT if they deviate from it due to disturbances or insufficient time allowed in the static start-up phase, 
{\vspace{1ex} \noindent \fbox{ \begin{minipage}{.95\columnwidth} 
\noindent
{\it TVD-C:}
\begin{equation}
\dot p = \left( I-\frac{\partial c }{ \partial p } \right)^{-1}\left(-k(p-c)+\frac{\partial c }{ \partial t } \right). \label{practical}
\end{equation}
\end{minipage}}\\[2ex]}
If a CVT is perfectly achieved initially, then the proportional term does not contribute anything to the update law, and the result in Theorem \ref{thm:thm1} still applies. We denote this algorithm {\it TVD-C}, where {\it TVD} stands for Time-Varying Densities, and {\it C} stands for Centralized, as will be discussed in subsequent sections.

The second issue with Equation \ref{eq:inv} is the presence of the term $\partial c/\partial p$. Even though this might look innocent, this term is rather complicated, due to the fact that
$$
c_i (p,t) = \frac{ \int_{V_i(p)} q \phi(q,t)dq}{\int_{V_i(p)} \phi(q,t)dq},
$$
which depends on $p$ in the boundary of the area over which the two integrals are taken. As a result, Leibniz rule must be exercised, and this computation is discussed in detail in Appendix A.

\section{DISTRIBUTED APPROXIMATIONS} \label{sec:Decentralization}
Given a Voronoi partition, we will denote the boundary between two cells by $\partial V_{ij}$. In the planar case, this boundary is either empty (Voronoi cells do not intersect), a single point (Voronoi cells intersect at a single vertex) or a line (Voronoi cells share a face). The two  Voronoi cells are said to be adjacent if they share a face, and we denote the set of cells adjacent to cell $i$ by $N_{V_i}$. 

Now, suppose that $i \not \in N_{V_j}$. This means either $\partial V_{i,j}$ is empty or consists of a singleton. This moreover implies that any integrals over $\partial V_{i,j}$ are zero, and (see Appendix A), Leibniz Rule tells us that these integrals are what define $\frac{\partial c_i }{ \partial p_j }$, from which we can conclude that $\frac{\partial c_i }{ \partial p_j } = 0$. As such we have the following result
\begin{lemma}
$i \not \in N_{V_j}\Rightarrow\frac{\partial c_i }{ \partial p_j } = 0$.
\label{lem:sparse}
\end{lemma}

A direct consequence of Lemma \ref{lem:sparse} is that $\partial c/\partial p$ encodes adjacency information. And, both algorithms in Equation \ref{eq:Lloyd} and \ref{eq:Cortes} are distributed in this manner, i.e., the update rule for $\dot p_i$ only depends on $p_j$ if $j\in N_{V_i}$. This, however, is not the case with Equation \ref{practical} since even though $\partial c/\partial p$ has the right sparsity structure, $(I-\partial c/\partial p)^{-1}$ does not. In fact, the inverse renders the resulting matrix dense and all sparsity structure is lost. The purpose of this section is thus twofold: 1) to develop a distributed approximation to \ref{practical} and 2) to overcome the restrictrions associated with $\phi$ for the inverse in \ref{practical} to exist. The required approximation can be found in the Neumann series, e.g., \cite{stewart1998matrix}.

\begin{lemma}[Neumann series]
Let A be a square matrix. If $\lim_{k \rightarrow \infty} A^k = 0$, then $I-A$ is invertible and 
\[
(I-A)^{-1} = I + A + A^2 + A^3 + \ldots.
\]
\end{lemma}
Moreover, for a $m \times m$ square matrix $A$, $\lim_{k \rightarrow \infty} A^k = 0$ if and only if $\abs{\lambda_i} < 1$ for all $i = 1, 2, \cdots, m$, where $\lambda_i$ are the eigenvalues of $A$.  As such, let  $\lambda_{max}$ denote the eigenvalue with the largest magnitude of the matrix ${\partial c}/{\partial p}$. Using the Neumann series, we can express $(I-{\partial c}/{\partial p})^{-1}$ as
\[
\br{I-\f{\partial c}{\partial p}}^{-1} = I + \f{\partial c}{\partial p} + \br{\f{\partial c}{\partial p}}^2 + \ldots
\]
as long as $\abs{ \lambda_{max} } < 1$. 

Now, if we insist on only letting $\dot p_i$ depend on $p_j,~j\in N_{V_i}$, (as well as $p_i$ itself) we can truncate the series after just two entries
\[
\br{I-\f{\partial c}{\partial p}}^{-1} \approx I + \f{\partial c}{\partial p},
\]
which gives the update law (modified from \ref{practical}),
	\[
	\dot p = \br{ I + \f{\partial c}{\partial p} }\left(-k(p-c) +\f{\p c}{\p t}\right),
	\] 
or at the level of the individual robots

{\vspace{1ex} \noindent \fbox{ \begin{minipage}{.95\columnwidth} 
\noindent
{\it TVD-}$D_1:$
\begin{equation}
\dot p_i=\f{\partial c_i}{\partial t} - k(p_i-c_i)+ \sum_{j\in N_{V_i}} \f{\partial c_i}{\partial p_j} \left(\f{\partial c_j}{\partial t}-k(p_j-c_j)\right),
\label{decentralized}
\end{equation}
\end{minipage}}\\[2ex]}
where the label denotes Time-Varying-Density, Decentralized with 1-hop adjacency information.

It should be noted that Equation \ref{decentralized} is always well-defined (as long as $\phi$ is continuously differentiable). In other words, even if the Neumann series is not convergent or if the inverse does not exist, the entries in \ref{decentralized} are well-defined. And, as will be seen in subsequent sections, even if $|\lambda_{max}|>1$ for short time-periods, it never stays that way as the robots evolve, meaning that the approximation is sound, eventually.

One can now investigate what happens when higher order terms are kept in the Neumann series. For this, we let $dist(i,j)$ denote the distance between cells $i$ and $j$. \begin{footnote}{Formally speaking, $dist(i,j)$ is the edge distance between $i$ and $j$ in the Delaunay graph induced by the Voronoi tessellation.}\end{footnote} And, as $\partial c/\partial p$ is a (block) adjacency matrix, 
we have that
$$
\left[\left(\partial c/\partial p\right)^k\right]_{ij}\neq0~\Rightarrow~dist(i,j)=k,~k=0,1,2,\ldots
$$
where $[\cdot]_{ij}$ denotes the block corresponding to cell $c_i$ and robot position $p_j$.

The $k$-hop version of $TVD-D_1$ thus becomes

{\vspace{1ex} \noindent \fbox{ \begin{minipage}{.95\columnwidth} 
\noindent
{\it TVD-}$D_k:$
\begin{equation}
\dot p = \sum_{\ell=0}^k\left(\frac{\partial c}{\partial p}\right)^\ell\left(-k(p-c) +\f{\p c}{\p t}\right),
\label{eq:khop}
\end{equation}
\end{minipage}}\\[2ex]}

In the next section, we implement these algorithms on a team of mobile robots, as well as discuss the implementation details
and a comparison to other methods.

\section{IMPLEMENTATION} \label{sec:rob}
In the previous section, a family of distributed algorithms, $TVD$-$D_k,~k=0,1,\ldots$, where developed as approximations to $TVD$-$C$, which, in turn, was presented as an alternative to the two algorithms dubbed $Lloyd$ and $Cortes$. In this section, we implement these different algorithms on a team of mobile robots, both in simulation and on Khepera III differential-derive mobile robots. Of particular importance is the fact that the Neumann approximation is really only valid when $|\lambda_{max}|<1$. But, as will be seen empiricially, this holds almost all the time, and when it does not, the robots quickly move into a configuration where the Neumann approximation is indeed valid.

\begin{figure}[ht!]%
\begin{center}
	\subfloat[]{\includegraphics[width=6cm]{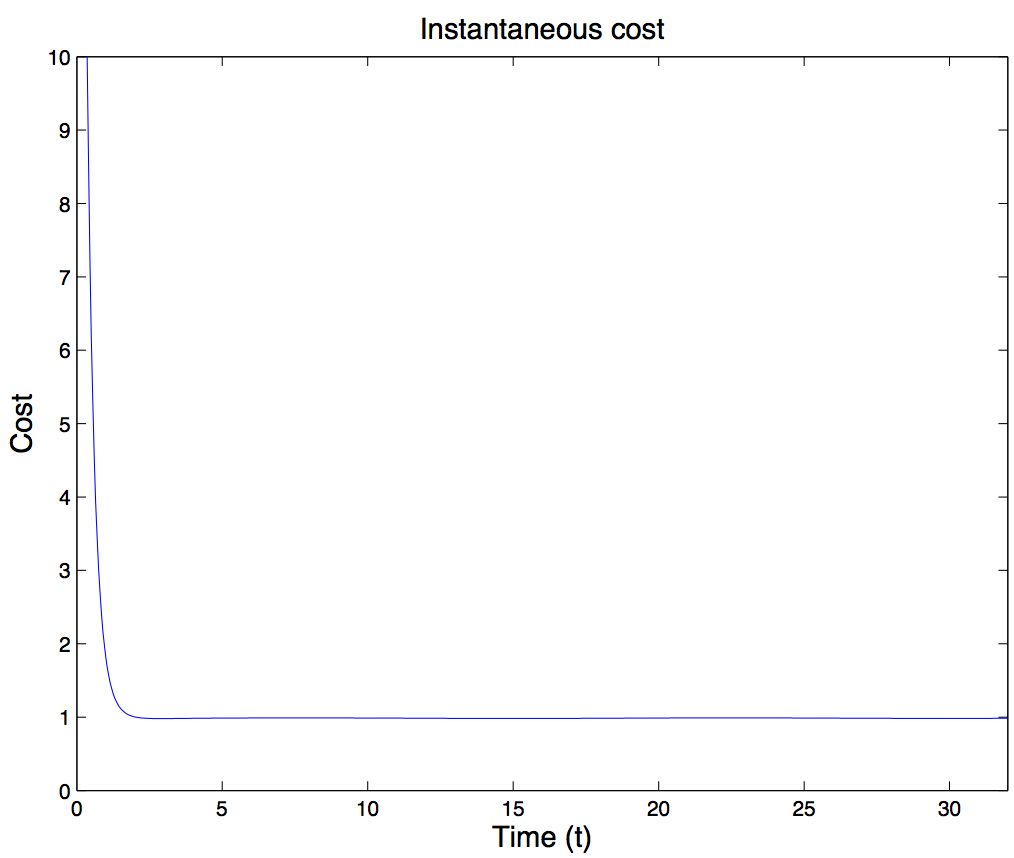}}\\
\subfloat[]{\includegraphics[width=6.1cm]{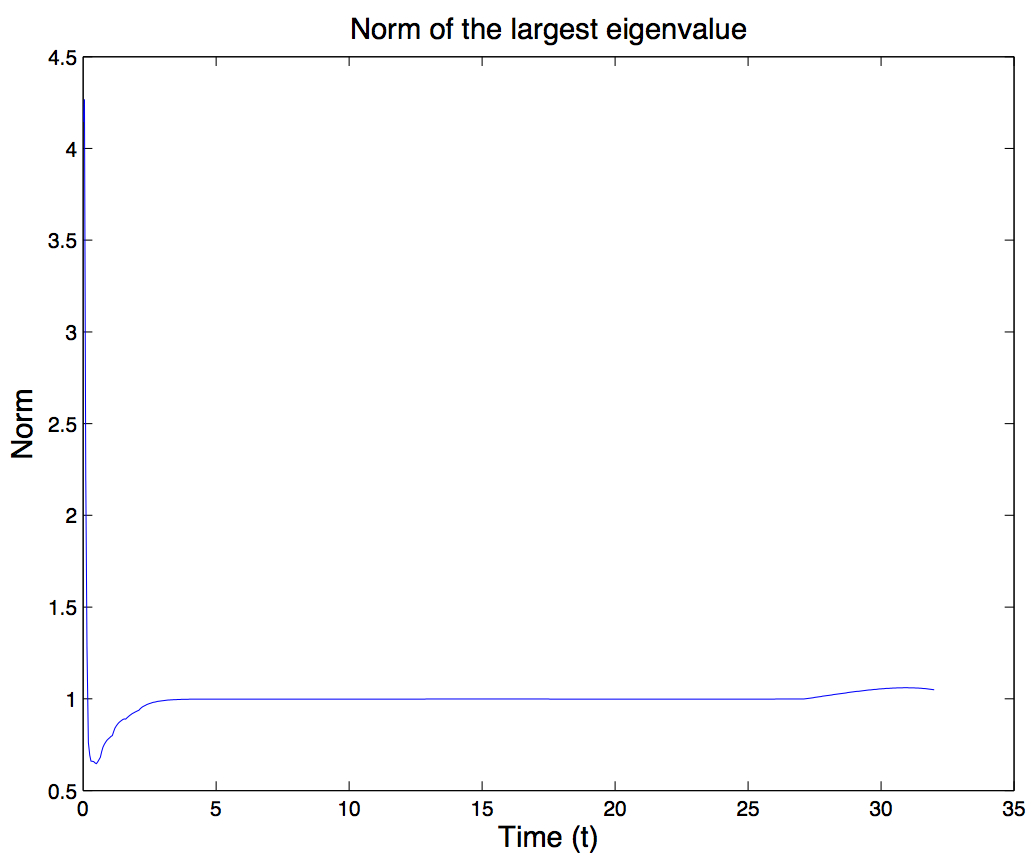}}
\end{center}
	\caption{Instantaneous locational cost (a) and the magnitude of $\lambda_{max}$ (b) as a function of time when executing $TVD$-$D_1$ under density function $\phi_2$.}
	\label{example1}
	\end{figure}

Two different density functions were considered
	\begin{align*}
	\phi_1(q,t) & = e^{ -\br {\br{ q_x-2\sin\br{\f{t}{\tau}} }^2  + \br{\f{q_y}{4}}^2 }}\\
\phi_2(q,t) & = e^{ -\br { \br{ q_x-2\cos\br{\f{t}{\tau}} }^2  + \br{ q_y-2\sin\br{\f{t}{\tau}} }^2 }} 
	\end{align*}
The time constant was taken to be $\tau=5$, and as a sanity-check, a number of simulations were performed using $TVD$-$D_1$ from different initial conditions, to gauge when the Neumann approximation was invalid and whether or not this affected the performance of the algorithm. As is shown in Figures \ref{example1}, this was not the case even when the magnitude of the maximal eigenvalue hovers close to 1 for longer periods of time.


Moreover, different verions of $TVD$-$D_k$, where simulated for $\phi_1$ and $\phi_2$ with the total cost being
$$
\int_0^{T_f}H(p(t),t)dt.
$$
The costs are summarized in Table \ref{table:neumann}. And, as can be seen, the cost does indeed decrease slightly as more terms are kept in the Neumann serier. However, the difference between the different cases is not particularly dramatic beyond the $k=0$ to $k=1$ case, i.e., when no information is used about neighboring robot positions and when only adjacent neighbors are taken into account. Similarly, the price of anarchy, i.e., the difference between $TVD$-$D_1$ and $TVD$-$C$ is marginal.

\begin{table}[!h]
\caption{Total costs under different $TVD$-$D_k$}
\begin{center}
    \begin{tabular}{| l | l | l |}
    \hline
Algorithm &	Total cost ($\phi_1$) & Total cost ($\phi_2$) \\ \hline
$TVD$-$D_0$ &	316.7 &	37.3\\ \hline
$TVD$-$D_1$ &	309.8 &	35.9\\ \hline
$TVD$-$D_2$ &	308.2 &	35.8\\ \hline
$TVD$-$D_3$ &	307.5 &	35.8\\ \hline
$TVD$-$D_4$ &	307.1 &	35.8\\ \hline
$TVD$-$D_5$ &	306.8 &	35.8\\ \hline
$TVD$-$D_6$ &	306.7 &	35.8\\ \hline
$TVD$-$D_7$ &	306.7 &	35.8\\ \hline
$TVD$-$D_8$ &	306.6 &	35.8\\ \hline
$TVD$-$D_9$ &	306.6 &	35.8\\ \hline
$TVD$-$D_{10}$ &	306.5 &	35.8\\ \hline
$TVD$-$C$  &	306.4 &	35.8\\ \hline
    \end{tabular} 
\end{center}
\label{table:neumann}
\end{table}

Moreover, a comparison was mande to $LLoyd$ and to $Cortes$, using $\phi_1$ and $\phi_2$ as well as three additional time-varying density functions. In all of these cases, the robots were initialized at the same positions to provide an inherently problematic comparison since the algorithms are chasing local (as opposed to global) minimizers to the locational cost.  These findings are summarized in Table \ref{table:algorithms}.

\begin{table}[!ht]
\caption{Coverage performance comparison.}
\begin{center}
    \begin{tabular}{| l | l | l | l | l | l |}
    \hline
	& $\phi_1$ &   $\phi_2$ & $\phi_3$ & $\phi_4$ &  $\phi_5$  \\ \hline
$TVD$-$D_{1}$	 	&	309.8   & 	35.0 &	36.5 &	35.9 &	100.2\\ \hline
$TVD$-$C$ &	306.4 &	34.3 &	34.3 &	36.3 &	98.9\\ \hline
${Cortes}$ 				&	319.5 &	38.4 &	N/A	    &	37.5 &	101.7\\ \hline
${Lloyd}$					&	324.6 &	40.1 &	52.6 &	38.7 &	103.6\\ \hline
    \end{tabular} 
\end{center}
\label{table:algorithms}
\end{table}

In all cases (except $\phi_3$) $Cortes$ did indeed perform better than $Lloyd$, which is not surprising since $Lloyd$ is designed for static density functions. However, under density function $\phi_3$, the assumptions behind $Cortes$ were violated and, as a result, the robots ended up leaving the domain $D$ when running $Cortes$. Moreover, $TVD$-$C$ and $TVD$-$D_1$ both outperformed $Cortes$ and $Lloyd$ in all five cases. In all of those cases, $TVD$-$C$ performed best, as can be expected, except for $\phi_4$, when $TVD$-$D_1$ was most effective. The reason for this is that the inverse was momentarily ill-defined in that particular case, causing the robots to move rather erratically during a short period of time when executing $TVD$-$C$. The outcome of this is that we prescribe $TVD$-$D_1$ as the overall most effective algorithm
 since it is always well-posed, allows for distributed implementation, and performs better than previously proposed algorithms.

\begin{figure*}%
\centering{
\subfloat[]{\includegraphics[width=8.0cm]{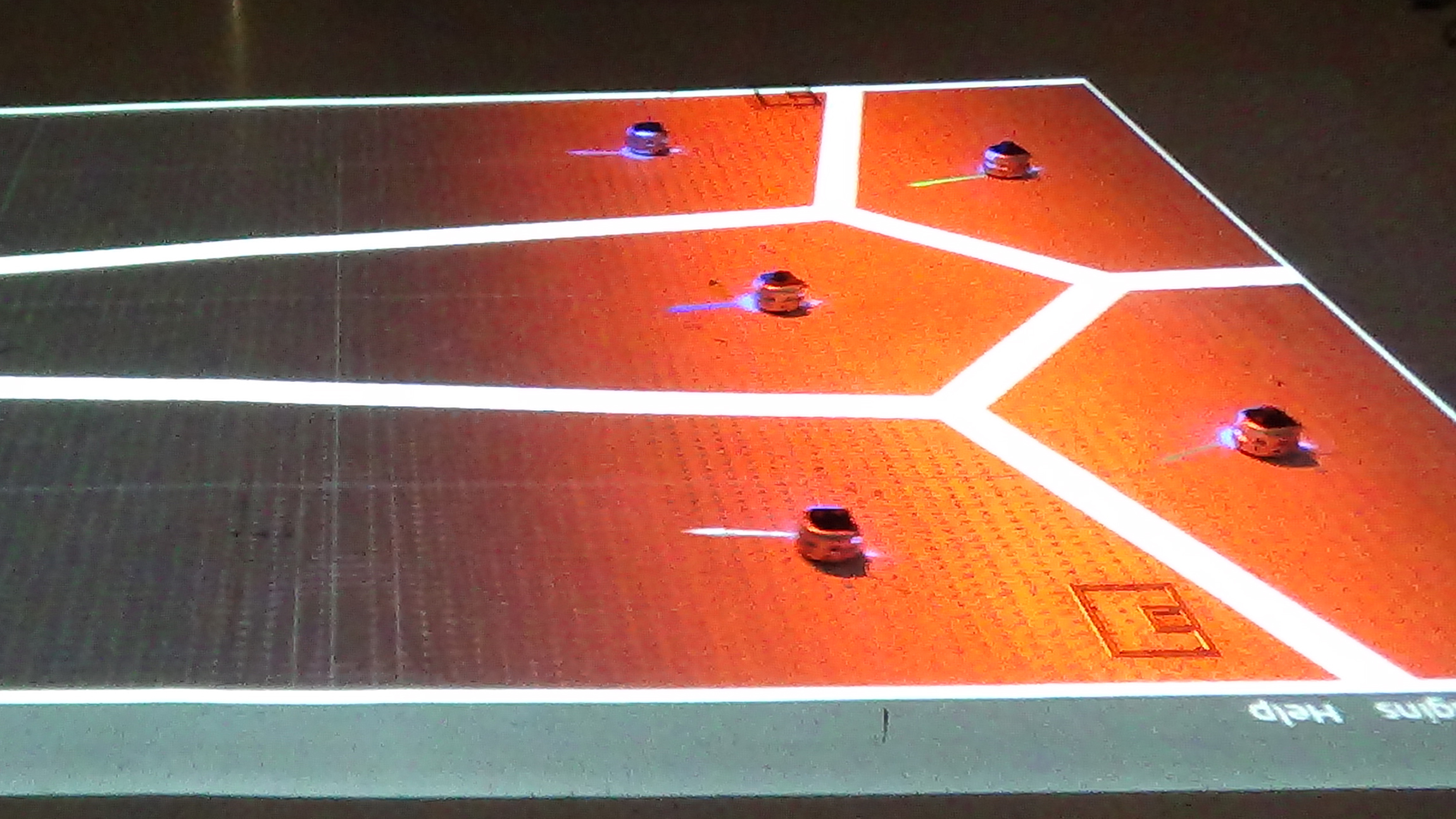}}\,
\subfloat[]{\includegraphics[width=8.0cm]{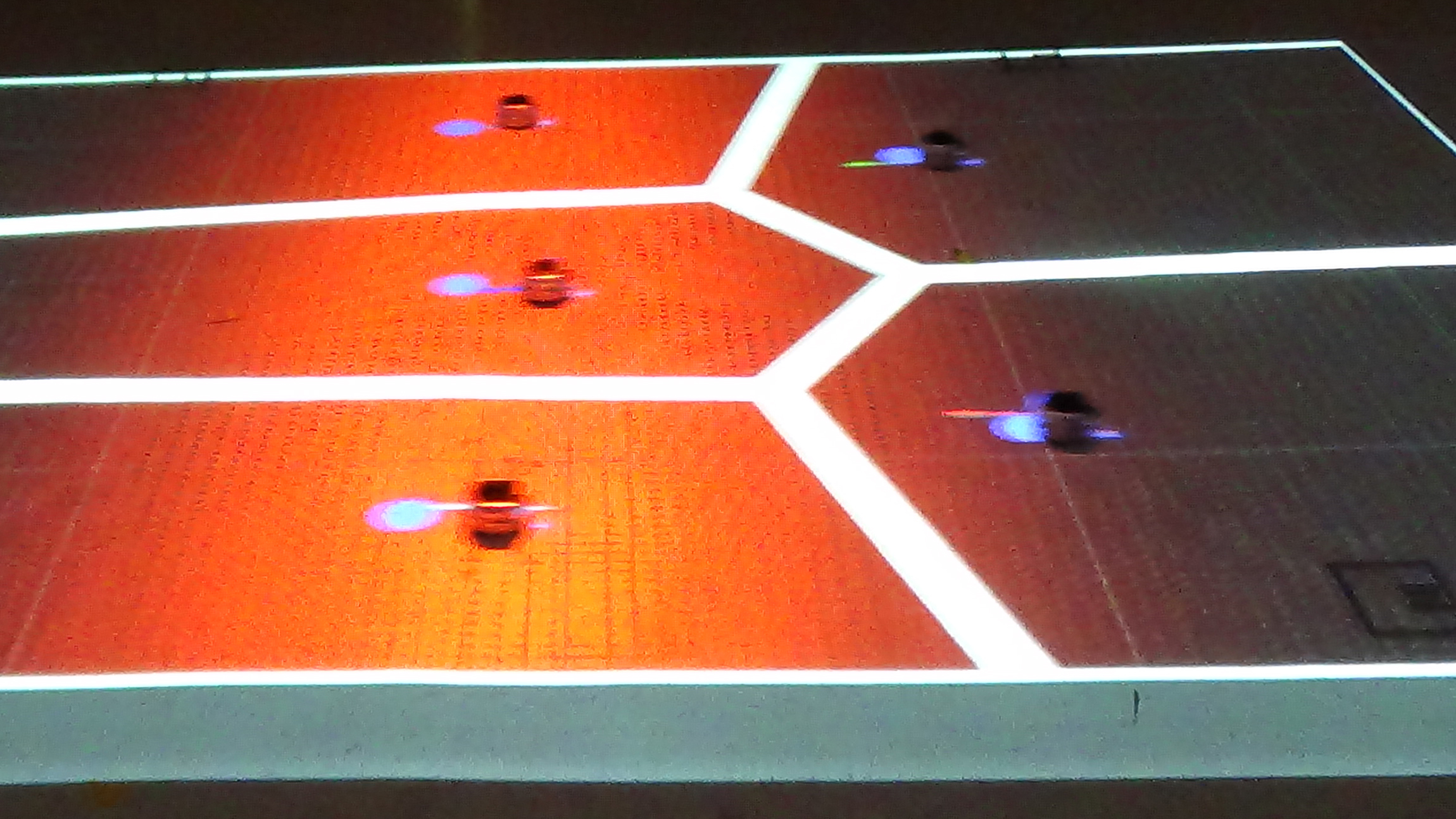}}\qquad
\subfloat[]{\includegraphics[width=8.0cm]{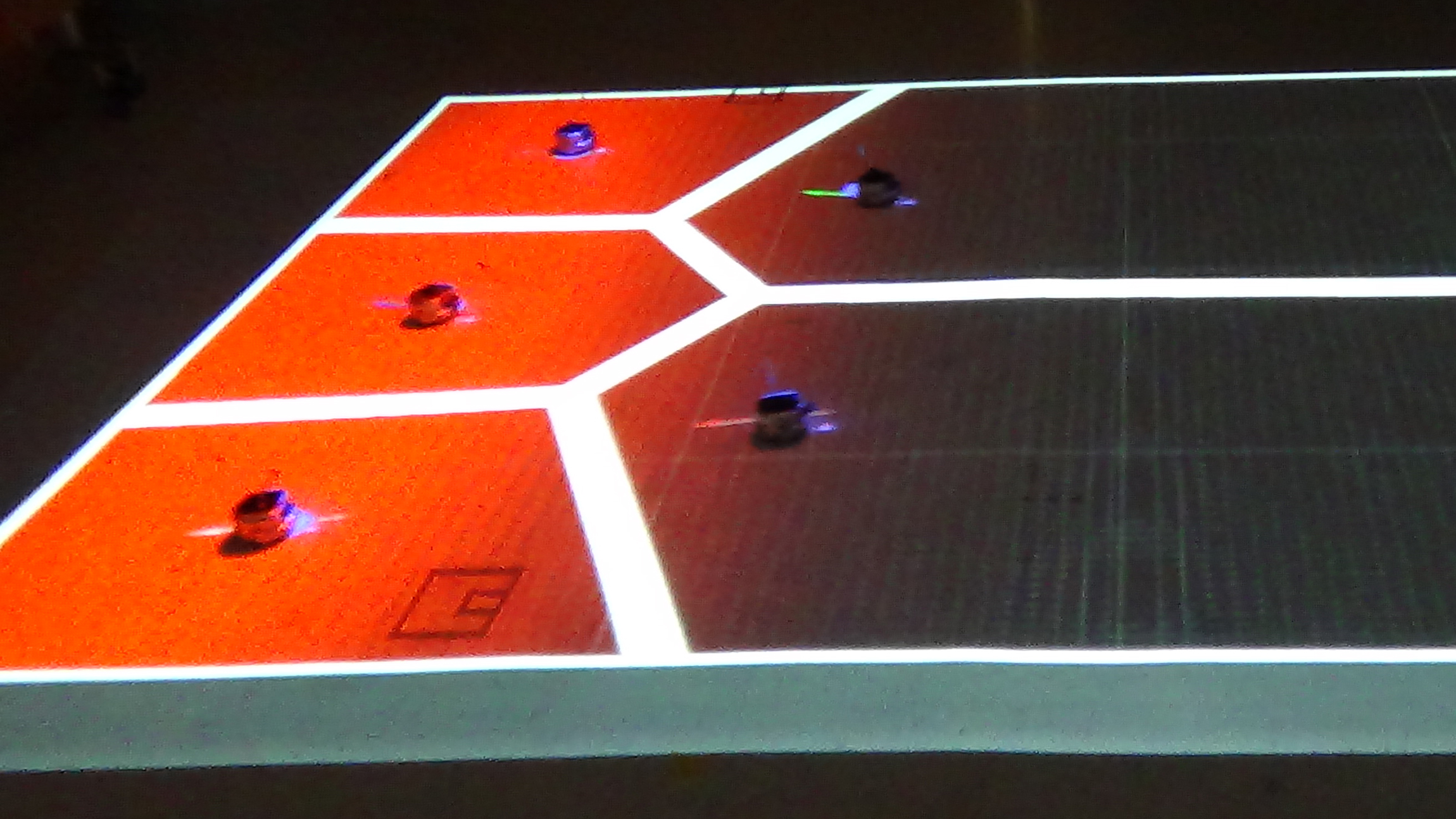}}\,
\subfloat[]{\includegraphics[width=8.0cm]{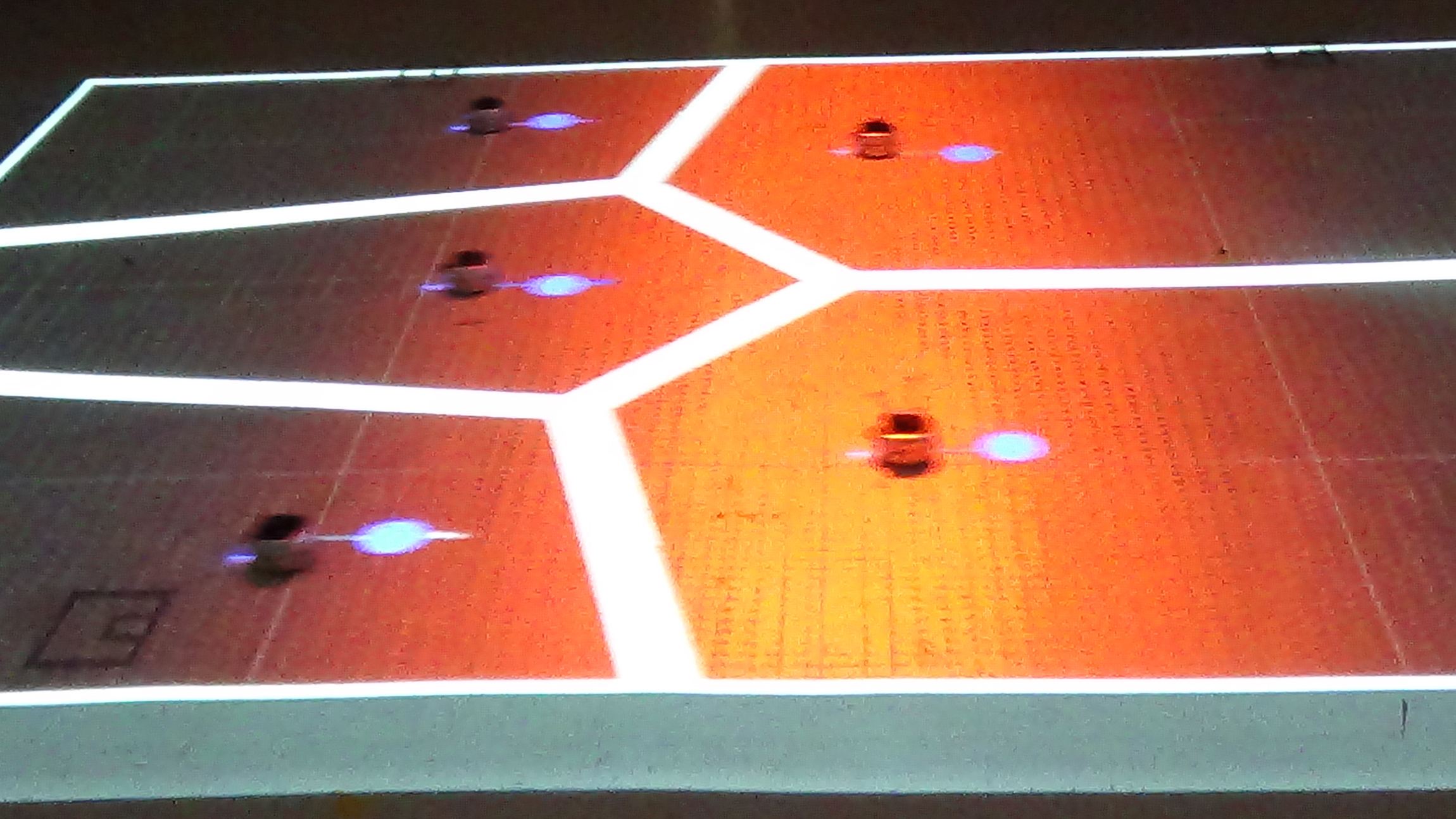}}\qquad
\subfloat[]{\includegraphics[width=8.0cm]{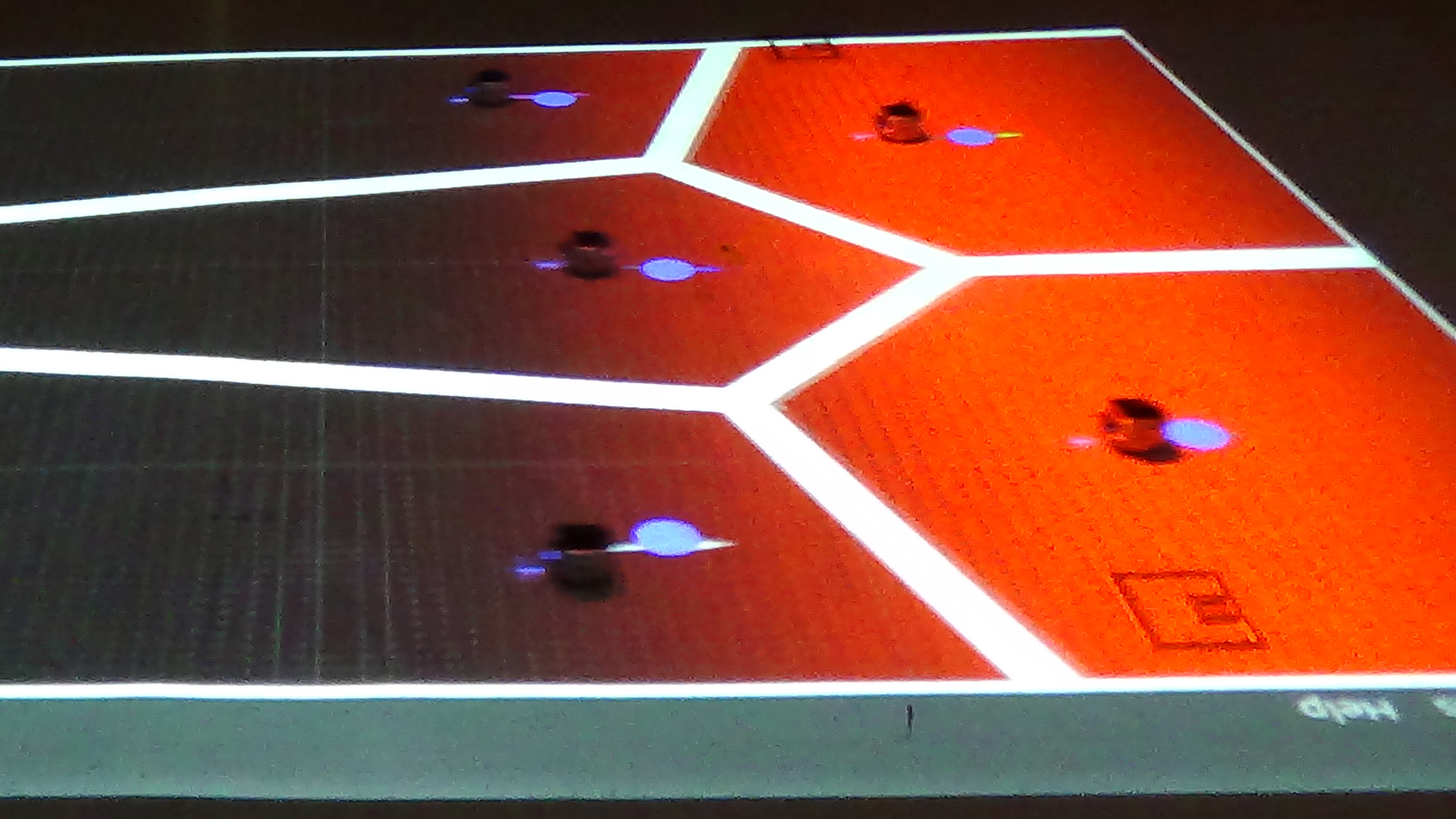}}\,
\subfloat[]{\includegraphics[width=8.0cm]{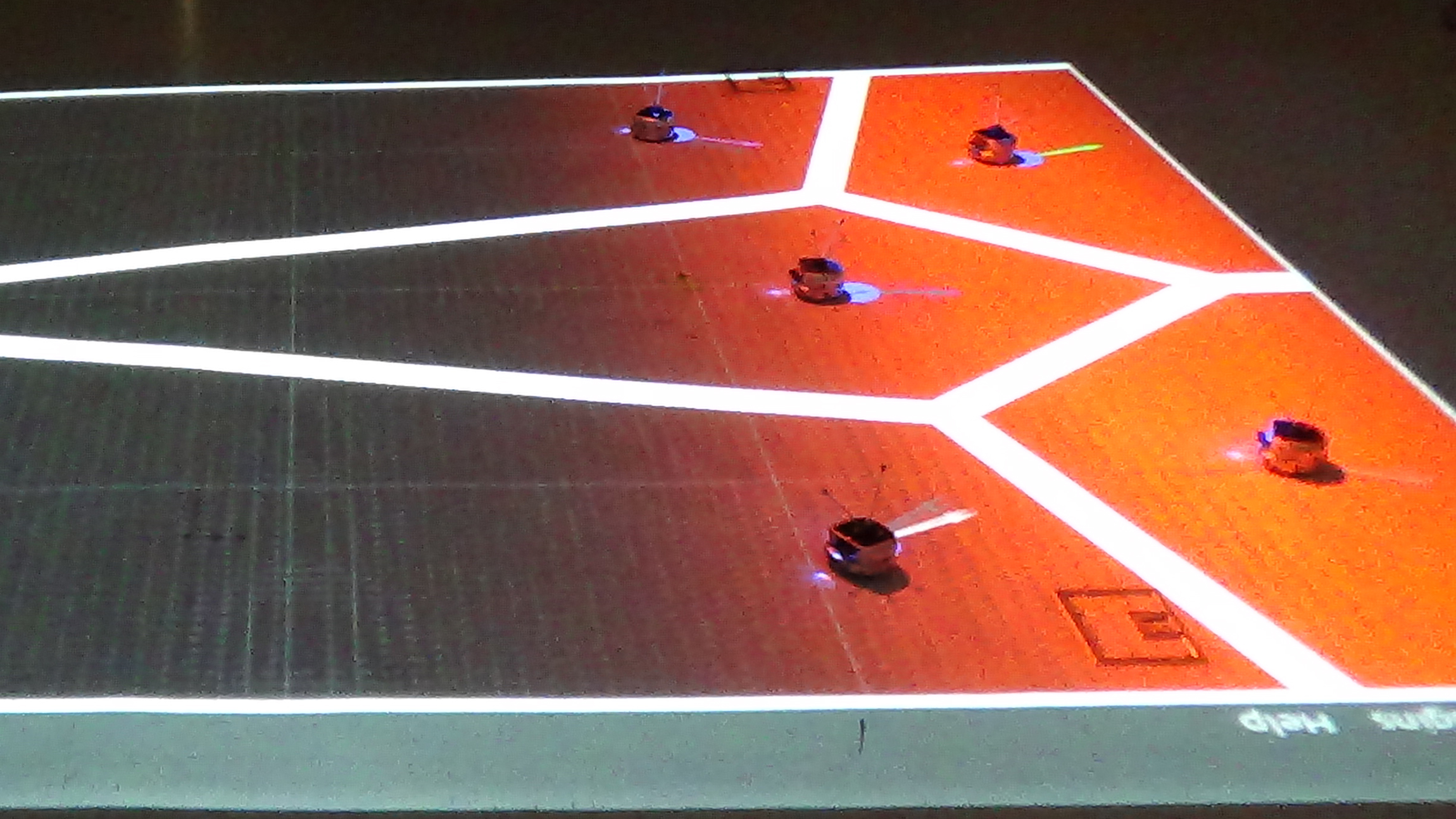}}\qquad}
\caption{Distributed coverage algorithm $TVD$-$D_1$ is deployed on a team of five mobile robots. An overhead projector is visualizing pertinent information, where the thick lines delineate the Voronoi cells, whose centers of mass are shown as bright circles. The corresponding video is available at
\texttt{https://www.youtube.com/watch?v=80YAsC3wVIk}}
\label{rviz}
\end{figure*}

As $TVD$-$D_1$ was the the all-around most effective algorithm based on the simulation results, this algorithm was implemented on a team of mobile robots. The ROS (Robot Operating System, version Diamondback) framework running on Ubuntu (version 11.04) machine with Intel dual core CPU 2.13GHz, 4GB memory was used to implement the algorithm and send control signals to individual robots over a wireless router. Five Khepera III robots from K-team were used as the team of mobile robots for the experiment. The Khepera III robots each have a 600MHz ARM processor with 128Mb RAM, embedded Linux, differential drive wheels, and a wireless card for communication over a wireless router. Ten Optitrack S250e motion capture cameras were used to provide position  and orientation data for the robots, which were used to provide the information required for the algorithm and the computation of the Voronoi partitions. The rviz package in ROS was used for visualizations, such as the position and the orientation of the robots, the density function, and the Voronoi partitions. The visualization was overlapped with the real physical environment to give a real-time visual representation,  as shown in Figure \ref{rviz}.

As the Khepera III mobile robots are differential-drive robots, they can be modeled as unicycles,
\begin{align*}
&\dot x_i = v_i\cos\theta_i\\
&\dot y_i = v_i\sin\theta_i\\
&\dot \theta_i = \omega_i,
\end{align*}
where $(x_i,y_i)$ is the position of robot $i$, $\theta_i$ its heading, and $v_i, \omega_i$ are the translational and angular velocities. In contrast to this, the coverage algorithm provides desired motions in terms of $\dot p_i$ and we map these onto $v_i,\omega_i$ through
\begin{align*}
&v_i = \norm{\dot p_i}\\
&\omega_i = \begin{bmatrix}-\sin\theta_i&\cos\theta_i\end{bmatrix} \cdot \f{\dot p_i}{\norm{\dot p_i}}.
\end{align*}

\section{CONCLUSIONS}

In this paper, we design coverage algorithms that allow for time-varying density functions. This is motivated by a desire to enable human operators to influence large teams of mobile robots, and density functions would constitute an abstraction for this that does not depend on the team size. Two different algorithms were considered, a centralized one that assumes that the agents already start at a centroidal Voronoi tessellation, and a distributed algorithm that does not require this constraint. In simulation, it is shown that the proposed methods outperform previously proposed methods, and the decentralized algorithm is deployed on a team of mobile robots.

\bibliographystyle{IEEEtran}
\bibliography{bibfile}

\begin{thebibliography}{10}
\providecommand{\url}[1]{#1}
\csname url@samestyle\endcsname
\providecommand{\newblock}{\relax}
\providecommand{\bibinfo}[2]{#2}
\providecommand{\BIBentrySTDinterwordspacing}{\spaceskip=0pt\relax}
\providecommand{\BIBentryALTinterwordstretchfactor}{4}
\providecommand{\BIBentryALTinterwordspacing}{\spaceskip=\fontdimen2\font plus
\BIBentryALTinterwordstretchfactor\fontdimen3\font minus
  \fontdimen4\font\relax}
\providecommand{\BIBforeignlanguage}[2]{{%
\expandafter\ifx\csname l@#1\endcsname\relax
\typeout{** WARNING: IEEEtran.bst: No hyphenation pattern has been}%
\typeout{** loaded for the language `#1'. Using the pattern for}%
\typeout{** the default language instead.}%
\else
\language=\csname l@#1\endcsname
\fi
#2}}
\providecommand{\BIBdecl}{\relax}
\BIBdecl

\bibitem{Kolling:HSI}
A.~Kolling, K.~Sycara, S.~Nunnally, and M.~Lewis, ``Human swarm interaction: An
  experimental study of two types of interaction with foraging swarms,''
  \emph{Journal of Human-Robot Interaction}, vol.~2, no.~2, pp. 103--128, June
  2013.

\bibitem{Martinez:Motion}
S.~Martinez, J.~Cortes, and F.~Bullo, ``{Motion Coordination with Distributed
  Information},'' \emph{Control Systems Magazine, IEEE}, vol.~27, no.~4, pp.
  75--88, 2007.

\bibitem{Ghosh:Coverage}
A.~Ghosh and S.~K. Das, ``Coverage and connectivity issues in wireless sensor
  networks: A survey.'' \emph{Pervasive and Mobile Computing}, vol.~4, no.~3,
  pp. 303--334, 2008.

\bibitem{Cortes:Coverage}
J.~Cortes, S.~Martinez, T.~Karatas, and F.~Bullo, ``Coverage control for mobile
  sensing networks: Variations on a theme,'' in \emph{Mediterranean Conference
  on Control and Automation}, Lisbon, Portugal, 2002, {Electronic Proceedings}.

\bibitem{Cortes:Coverage2}
------, ``Coverage control for mobile sensing networks,'' \emph{IEEE
  Transactions on Robotics and Automation}, vol.~20, no.~2, pp. 243--255, Apr.
  2004.

\bibitem{Cortes:Coordination}
J.~Cort{\'e}s and F.~Bullo, ``Coordination and geometric optimization via
  distributed dynamical systems,'' \emph{SIAM Journal on Control and
  Optimization}, vol.~44, no.~5, pp. 1543--1574, October 2005.

\bibitem{Pimenta:SCAT}
L.~C.~A. Pimenta, M.~Schwager, Q.~Lindsey, V.~Kumar, D.~Rus, R.~C. Mesquita,
  and G.~A.~S. Pereira, ``Simultaneous coverage and tracking (scat) of moving
  targets with robot networks.'' in \emph{WAFR}, ser. Springer Tracts in
  Advanced Robotics, vol.~57.\hskip 1em plus 0.5em minus 0.4em\relax Springer,
  2008, pp. 85--99.

\bibitem{Haumann:Discoverage}
D.~Haumann, V.~Willert, and K.~D. Listmann, ``Discoverage: From coverage to
  distributed multi-robot exploration,'' in \emph{Proceedings of the 4th IFAC
  Workshop on Distributed Estimation and Control in Networked Systems}, vol.~4,
  no.~1, Koblenz, Germany, September 2013, pp. 328--335.

\bibitem{Macwan:Target}
A.~Macwan, G.~Nejat, and B.~Benhabib, ``Target-motion prediction for robotic
  search and rescue in wilderness environments.'' \emph{IEEE Transactions on
  Systems, Man, and Cybernetics, Part B}, vol.~41, no.~5, pp. 1287--1298, 2011.

\bibitem{Ghaffarkhah:DynamicCoverage}
A.~Ghaffarkhah, Y.~Yan, and Y.~Mostofi, ``Dynamic coverage of time-varying
  environments using a mobile robot -- a communication-aware perspective,'' in
  \emph{GLOBECOM Workshops (GC Wkshps), 2011 IEEE}, 2011, pp. 1297--1302.

\bibitem{Du:CVT}
Q.~Du, V.~Faber, and M.~Gunzburger, ``Centroidal voronoi tessellations:
  Applications and algorithms,'' \emph{SIAM Review}, vol.~41, no.~4, pp.
  637--676, Dec. 1999.

\bibitem{Meguerdichian:Exposure}
S.~Meguerdichian, F.~Koushanfar, G.~Qu, and M.~Potkonjak, ``Exposure in
  wireless ad-hoc sensor networks,'' in \emph{Proceedings of the 7th annual
  international conference on Mobile computing and networking}, ser. MobiCom
  '01.\hskip 1em plus 0.5em minus 0.4em\relax New York, NY, USA: ACM, 2001, pp.
  139--150.

\bibitem{Adlakha:Critical}
S.~Adlakha and M.~B. Srivastava, ``Critical density thresholds for coverage in
  wireless sensor networks.'' in \emph{WCNC}.\hskip 1em plus 0.5em minus
  0.4em\relax IEEE, 2003, pp. 1615--1620.

\bibitem{Iri:Fast}
M.~Iri, K.~Murota, and T.~Ohya, ``A fast voronoi-diagram algorithm with
  applications to geographical optimization problems,'' in \emph{System
  Modelling and Optimization}, ser. Lecture Notes in Control and Information
  Sciences, P.~Thoft-Christensen, Ed.\hskip 1em plus 0.5em minus 0.4em\relax
  Springer Berlin Heidelberg, 1984, vol.~59, pp. 273--288.

\bibitem{Du:Convergence}
Q.~Du, M.~Emelianenko, and L.~Ju, ``Convergence of the lloyd algorithm for
  computing centroidal voronoi tessellations,'' \emph{SIAM Journal on Numerical
  Analysis}, vol.~44, no.~1, pp. 102--119, Jan. 2006.

\bibitem{Liu:CVT}
Y.~Liu, W.~Wang, B.~L{\'e}vy, F.~Sun, D.-M. Yan, L.~Lu, and C.~Yang, ``On
  centroidal voronoi tessellation -- energy smoothness and fast computation,''
  \emph{ACM Transactions on Graphics}, vol.~28, no.~4, pp. 1--17, Sep. 2009.

\bibitem{Du:Acceleration}
Q.~Du and M.~Emelianenko, ``Acceleration schemes for computing centroidal
  voronoi tessellations,'' \emph{Numerical Linear Algebra with Applications},
  vol.~13, no. 2-3, pp. 173--192, 2006.

\bibitem{stewart1998matrix}
G.~Stewart, \emph{Matrix Algorithms Volume 1: Basic Decompositions}.\hskip 1em
  plus 0.5em minus 0.4em\relax Society for Industrial and Applied Mathematics,
  1998.

\end{thebibliography}

\section*{APPENDIX A}

Central to the developments in this paper is the computation of the partial derivative $\frac{\partial c }{ \partial p }$. The first step towards this is the application of Leibniz rule, e.g., \cite{Du:Acceleration}.
\begin{lemma}\label{partial derivative lemma1}
\textit{
Let $\Omega(p)$ be a region that is a smooth function of $p$ such that the unit outward normal vector $n$ is uniquely defined almost everywhere on $\partial \Omega$, which is the boundary of $\Omega$. Let
\[
F = \int_{\Omega(p)} f(q)dq.
\]
Then 
\[
\frac{\partial F}{\partial p} = \int_{\partial \Omega(p)} f(q)\hat q \cdot n(q) dq
\]
where $\hat q$ is the derivative of the points on $\partial \Omega$ with respect to $p$.
}
\end{lemma}

In \cite{Du:Acceleration}, it was investigated how Voronoi cells changed as functions of $p_i$. In fact, it was shown in \cite{Du:Acceleration} that
for any point $q \in \p V_{i,j}$ (the boundary between adjacent cells $i$ and $j$),
\begin{align*}
\frac{\partial q}{\partial p_j^{(b)}} \cdot (p_j - p_i) = \frac{1}{2} e_b\cdot(p_j - p_i) - e_b\cdot\left(q-\frac{p_i+p_j}{2}\right),\\
\frac{\partial q}{\partial p_i^{(b)}} \cdot (p_j - p_i) = \frac{1}{2} e_b\cdot(p_j - p_i) + e_b\cdot\left(q-\frac{p_i+p_j}{2}\right),
\end{align*}
where $p_j^{(b)}$ denotes the $b$-th component of the vector $p_j$ and $e_b$ is the $b$-th elementary unit vector, with $b=1,2$ in the planar case (which is what is considered in this paper). Note that in this paper, $b = 1, 2$ since we are considering the case $D \subset \R^2$ only.

Substituting this into Leibniz rule, we obtain
\begin{equation}
\begin{split}
&\frac{\partial c_i^{(a)} }{ \partial p_j^{(b)} } = \left(  \int_{\p V_{i,j}} \phi q^{(a)}\frac{ p_j^{(b)}  - q^{(b)} }{\norm{p_j - p_i }} dq \right)  \bigg/ m_i  \\
& -\left( \int_{\p V_{i,j}} \phi \frac{ p_j^{(b)}  - q^{(b)} }{\norm{p_j - p_i }}  dq  \right)    \left( \int_{V_i (P)} \phi q^{(a)}dq  \right) \bigg/ {m_i^2}, \label{i-j}
\end{split}
\end{equation}
where $a=1,2$ and where $i\neq j$.
Similarly, when $i=j$ we get
\begin{equation*}
\begin{split}
&\frac{\partial c_i^{(a)} }{ \partial p_i^{(b)} } = \left(  \int_{\p V_{i,j}} \phi q^{(a)}\frac{ q^{(b)} - p_i^{(b)} }{\norm{p_j - p_i }}  dq \right)  \bigg/ m_i\\
& -\left( \int_{\p V_{i,j}} \phi  \frac{ q^{(b)} - p_i^{(b)} }{\norm{p_j - p_i }}  dq  \right)    \left( \int_{V_i (P)} \phi q^{(a)}dq  \right) \bigg/ {m_i^2},
\end{split}
\end{equation*}
which gives us all we need to compute $\partial c/\partial p$.

\end{document}